\documentclass{article}
\usepackage{graphicx}
\usepackage{amsmath}
\usepackage{amsthm}
\usepackage{caption}
\providecommand{\keywords}[1]{\textbf{Keywords} #1}
\newtheorem{remark}{Remark}

\newtheorem{defn}{Definition}

\begin{document}

\title{On subgrid multiscale stabilized finite element method for advection-diffusion-reaction equation with variable coefficients}

\author{ Manisha Chowdhury  ,  B.V. Rathish Kumar \thanks{ Corresponding author,\newline Email addresses: chowdhurymanisha8@gmail.com (M. Chowdhury) and drbvrk11@gmail.com (B.V.R. Kumar)  } }
      
\date{Indian Institute of Technology Kanpur \\ Kanpur, Uttar Pradesh, India}

\maketitle
\begin{abstract}
In this study a stabilized finite element method for solving advection-diffusion-reaction equation with spatially variable coefficients has been carried out. Here subgrid scale approach along with algebraic approximation to the sub-scales has been chosen as stabilized method among various other methods. Both a priori and a posteriori finite element error estimates in $L_2$ norms have been derived after introducing the stabilized formulation of the variational form. An expression of the stabilization parameter for this 2D problem has also been derived here. At last numerical experiments are presented to verify numerical performance of the stabilized method and check the credibility of the theoretically derived expression of the stabilization parameter.
\end{abstract}


\section{Introduction}
Due to wide range of applications of the transport equation in the fields of bio-medical engineering, chemical engineering, environmental sciences etc.,  solving more complex form of this equation has always been in great demand. One of the appropriate mathematical model equations expressing transport phenomena is the advection-diffusion-reaction (ADR) equation and its different forms. Here particularly we dealing with ADR equation with spatially variable diffusion coefficients and velocity components, due to its importance in modelling several physical phenomena in more effective way. Since finiding analytical solution of these type of equation is not an easy job, numerical methods have been preferred in solving it. The article \cite{RefE} has presented a priori and a posteriori error bounds for the numerical scheme, Galerkin finite element method. It is well-known that the Galerkin formulation suffers from lack of stability when the diffusion coefficients are small compared to either advection velocity components or reaction term. The purpose of this paper is to present a study of a stabilized finite element method for this particular type of ADR equation. Here we have considered to study the stabilization method named by algebraic subgrid scale (ASGS) method i.e. subgrid scale (SGS) method with an algebraic approximation to the subscales. \cite{RefA} gives an essence of stabilized multiscale subgrid method for Helmholtz equation, whereas \cite{RefB} presents SGS formulation of ADR equation with constant coefficients. We have followed the same considerations to have SGS formulation with spatially variable coefficients in section 2 of this paper. Before that we introduce Galerkin finite element formulation for the equation in the same section. In sections 3 and 4, a priori and a posteriori errors have been estimated for stabilized form respectively. We have taken the conception of estimating errors from \cite{RefG} which provides a tutorial on the derivation of a priori and explicit residual-based a posteriori error estimates for Galerkin finite element discretization of general linear elliptic operators and \cite{RefF} which presents residual-based a posteriori bound for Galerkin least-squares (GLS) formulation for Helmholtz equation. A subsection of the section 4 presents an important highlight of this paper. It elaborately establishes the condition on stabilization parameter for 2D ADR equation because study of stabilization parameter only for 1D ADR equation is available in literatures \cite{RefB},\cite{RefC}. In the last section numerical experiments have been presented to verify the theoretically established results of error estimations as well as to check the accuracy of the expression of the stabilization parameter for a test problem whose conception is taken from literature in Hydrology \cite{RefD}. This section also presents a comparison of stability between the Galerkin method and the SGS method for diffusion dominated and advection dominated flows and then numerically establishes the fact that the Galerkin method suffers from lack of stability.
\section{Statement of the problem}
Introducing the notations,\\
c: concentration of dispersing mass of the solute\\
$(D_1,D_2)$: diffusion coefficients along x-axis and y-axis\\
$(u_1,u_2)$: velocity components along x-axis and y-axis\\
$\mu$: coefficient of reaction\\
$q$: source of the solute mass in the flow domain\\ 

The two-dimensional ADR equation in general form with spatially variable coefficient functions representing transportation of a solute in a bounded domain $\Omega$ in $R^2$ along with homogeneous Dirichlet boundary condition can be written as follows:\\
Find c : $ \Omega \rightarrow $ R such that\\
\begin{equation}
     \frac{\partial}{\partial x}[D_1 \frac{\partial c}{\partial x}- u_1 c] + \frac{\partial}{\partial y}[D_2 \frac{\partial c}{\partial y} - u_2 c] - \mu c + q = 0 \hspace{2mm} in \hspace{1mm} \Omega
\end{equation}
\hspace{8.2cm} c = 0 on $\partial\Omega $ \vspace{3mm}\\
Let us introduce two notations, $\widetilde{\bigtriangledown}$ := $(D_1 \frac{\partial}{\partial x},D_2 \frac{\partial}{\partial y})$ and $\overline{u}=(u_1,u_2)$ \\
Considering the fluid as incompressible i.e. the divergence of advection velocity is zero, the above equation can be rewritten with the ADR operator $\mathcal{L}$,
\begin{equation}
\mathcal{L}  = -\bigtriangledown \cdot \widetilde{\bigtriangledown} + \overline{u} \cdot \bigtriangledown  + \mu  \\
\end{equation}
whose adjoint operator $\mathcal{L^*}$ is \\
\begin{equation}
\mathcal{L^*}  = -\bigtriangledown \cdot \widetilde{\bigtriangledown} - \overline{u} \cdot \bigtriangledown  + \mu 
\end{equation}
The weak form of the boundary value problem is as follows: Find c $\in$ V such that\\
\begin{equation}
  a(c,d)=l(d) \quad \forall d \in V  
\end{equation}
where \\

\hspace{0.7cm} a(c,d)=$\int_{\Omega}D_1\frac{\partial c}{\partial x}\frac{\partial d}{\partial x} + \int_{\Omega}D_2\frac{\partial c}{\partial y} \frac{\partial d}{\partial y} + \int_{\Omega}d $$\bar{u}$$ \cdot \nabla c  +\int_{\Omega} \mu c d  $\\ 

\hspace{1cm} l(d)=$\int_{\Omega} q d = (d,q) $ \vspace{3 mm}\\
and let V represent both the usual spaces of trial solutions and weight functions. In addition $(\cdot ,\cdot): V \times V \rightarrow R$ represents $L_2(\Omega)$ inner product and $\| \cdot\|$ denotes $L_2(\Omega)$ norm, introduced by this inner product, i.e.
$\|d\|^2 = (d,d)= \int_{\Omega} d^2 $.\\
The subscripts on the inner product and norm notations (whenever it will be used) denote the domain of integration other than $\Omega$.
\begin{remark}
$a(c,d)=(d,\mathcal{L}c)= (\mathcal{L^*}d,c)$
\end{remark}
\begin{remark}
Furthermore  a(d,c)=$\int_{\Omega}D_1\frac{\partial d}{\partial x}\frac{\partial c}{\partial x} + \int_{\Omega}D_2\frac{\partial d}{\partial y} \frac{\partial c}{\partial y} - \int_{\Omega}c $$\bar{u}$$ \cdot \nabla d  +\int_{\Omega} \mu c d  $\\
Hence the bilinear form a(c,d) is not symmetric i.e. a(c,d) $\neq$ a(d,c).\\

\end{remark}

\subsection{Galerkin finite element formulation}
Let the bounded domain $\Omega$ be discretized into finite number of element subdomains $\Omega_e$, for e=1,2,...,$n_{el}$, where $n_{el}$ is the number of elements. Let $h_e$ be the diameter of each element $\Omega_e$ and h =$\underset{e=1,2,...,n_{el}}{max}$ $h_e$ \vspace{2mm}\\
Let  $V_h$ $\subset$ V be a finite dimensional space containing continuous piecewise polynomials of order p.
Now the Galerkin formulation is : Find $c_h \in V_h$ such that\\
\begin{equation}
    a(c_h,d_h) = l(d_h) \quad \forall d_h \in V_h
\end{equation}
where \\

\hspace{0.7 cm} $a(c_h,d_h)$=$\int_{\Omega}D_1\frac{\partial c_h}{\partial x}\frac{\partial d_h}{\partial x} + \int_{\Omega}D_2\frac{\partial c_h}{\partial y} \frac{\partial d_h}{\partial y} + \int_{\Omega}d_h $$\bar{u}$$ \cdot \nabla c_h  +\int_{\Omega} \mu c_h d_h  $\\ 

\hspace{1cm} $l(d_h)$=$\int_{\Omega} q d_h $ \vspace{3 mm}\\
Let $n_{pt}$ be the total number of nodes after discretizing the domain and $S_i$ be the standard shape(basis) function at node each node a, for a=1,2,...,$n_{pt}$ . Therefore the functions in $V_h$ can be interpolated from these shape functions as,\vspace{2mm}\\
$c_h$= $\sum_{a=1}^{n_{pt}} c_a S_a$ \quad $d_h$ = $\sum_{b=1}^{n_{pt}} d_b S_b $\\
Therefore\\
\begin{equation}
    \begin{split}
        a(c_h,d_h) & = \sum_{a,b=1}^{n_{pt}} d_b \hspace{1 mm} a(S_a,S_b) \hspace{1mm} c_a = D^t A C \\
       l(d_h) & = \sum_{b=1}^{n_{pt}} d_b \hspace{1mm} l(S_b) = D^t L \\
    \end{split}
\end{equation}
where A is a matrix of order $n_{pt} \times n_{pt}$ with components $a(S_a,S_b)$ and L is a column vector with components $l(S_b)$ for a,b=1,2,...,$n_{pt}$\\
Hence the discrete problem is equivalent to solve as the following linear system
\begin{equation}
    AC=L
\end{equation}

\subsection{Stabilized integral form}
Let $\Omega'$ be the union of element interiors, that is $\Omega' = \cup_{e=1}^{n_{el}} \Omega_k$. The stabilized form of equation (1) after applying SGS method,following \cite{RefC}, on it, is stated as follows: Find $c_h \in V_h$ such that \\
\begin{equation}
a_{SGS} (c_h,d_h)=l_{SGS} (d_h) \quad \forall d_h \in V_h
\end{equation}
where \\

\hspace{0.7 cm} $a_{SGS} (c_h,d_h)$= $ a(c_h,d_h)+ \int_{\Omega'} (-\mathcal{L^*}d_h) \tau \mathcal{L}c_h  $\\ 

\hspace{1cm} $l_{SGS}(d_h)$= $l(d_h) + \int_{\Omega'} (-\mathcal{L^*}d_h) \tau q $ \vspace{3 mm}\\
where $\tau$ is known as the stabilization parameter, which is measured later in this paper.  

\section{A priori error estimation in the $L_2$ norm}
A priori estimate provides a glimpse of convergence of the approximation method though it cannot be computable as it depends upon the exact solution. It also gives order of convergence before the solution is known. Before going to estimate the error let us first introduce its splitting which plays an important role in estimation. Let $c_h$ be the computed finite element solution using stabilized SGS method and c be the exact solution. Due to discretization of the domain we have an opportunity to introduce nodal interpolants. Let $c_h'$ be the nodal interpolant of c. Let e denote the finite element error.\\
Now let us introduce the error splitting as follows:
\begin{equation}
\begin{split}
e &= c_h-c\\
  &= (c_h-c_h')+(c_h'-c)\\
  &= e_h + \lambda
  \end{split}
\end{equation} 
where the first part $e_h= (c_h-c_h')$ of the error belongs to $V_h$ and the other part $\lambda= (c_h'-c)$ is the interpolation error belonging to V.
\subsection{Nitsche trick}
Generally Nitsche trick starts with considering an auxiliary problem. Here this auxiliary problem involves the variational form of the ADR equation expressed by adjoint operator $\mathcal{L^*}$, considering error e as the source term. This is also known as dual problem. The reason behind replacing the source term q by the error e is to introduce $L_2$ norm of e , which will be seen in further proceedings. Now let us introduce the dual problem as follows: Find $\eta \in V$ such that \\
\begin{equation}
a(d,\eta)=(e,d) \quad \forall d \in V
\end{equation}
Now replacing d by e in the above equation, we will have
\begin{equation}
a(e,\eta)= (e,e)= \|e\|^2
\end{equation}
The exact solution will also satisfy both the Galerkin finite element formulation (5) and SGS formulation (8) and both will give the following result:
\begin{equation}
a(c,d_h)=l(d_h) \quad \forall d_h \in V_h
\end{equation}
Now subtracting (10) from (5), we will have 
\begin{equation}
a(e, d_h)=0 \quad \forall d_h \in V_h 
\end{equation}
This is also known as Galerkin orthogonality.\\
Now we are going to bring SGS formulation into estimation by subtracting (10) from (8) as follows:
\begin{equation}
\begin{split}
a(c_h-c,d_h)+ (-\mathcal{L^*} d_h,\tau \mathcal{L}c_h) & = (-\mathcal{L^*} d_h,\tau q) \quad \forall d_h \in V_h \\
a(e,d_h) + (-\mathcal{L^*} d_h,\tau \mathcal{L}c_h) & = (-\mathcal{L^*} d_h,\tau \mathcal{L}c ) \quad \forall d_h \in V_h\\
a(e,d_h)+ (-\mathcal{L^*} d_h,\tau \mathcal{L}e) & = 0 \quad \forall d_h \in V_h\\
a_{SGS}(e,d_h)& =0 \quad \forall d_h \in V_h \\
\end{split}
\end{equation}
Let $\eta_h$ be the interpolant of $\eta $ in $V_h $ and $\lambda_\eta$ = $(\eta_h- \eta)$  be the interpolation error.\\
Since $\eta_h$ is in $V_h$, we can rewrite the last equation as follows,
\begin{equation}
a(e,\eta_h)+ (-\mathcal{L^*} \eta_h,\tau \mathcal{L}e) = 0
\end{equation}
To estimate error in $L_2$ norm, now we are bringing it into the context by subtracting (15) from (11). Hence,
\begin{equation}
\begin{split}
\|e\|^2 &= a(e, \eta-\eta_h)- (-\mathcal{L^*} \eta_h,\tau \mathcal{L}e)\\
        &= (\mathcal{L^*} \eta_h,\tau \mathcal{L}e)- a(e, \lambda_\eta)\\
        &=(\mathcal{L^*} \eta_h,\tau \mathcal{L}e_h) + (\mathcal{L^*} \eta_h,\tau \mathcal{L}\lambda)- a(e, \lambda_\eta)\\
\end {split}
\end{equation}
We have to find bounds for each term of the above equation so that one $\|e\|$  will be cancelled out from both sides. To have such bounds we need few important results, which are given as follows:\\
Standard interpolation estimate \cite{RefG} is,
\begin{equation}
\|\lambda_c\|_s \leq \bar{C}(p,\Omega) h^{p+1-s} \|c\|_{p+1}
\end{equation}
The exact solution c has been assumed of regularity (p+1). In addition $\|\cdot\|_r$ denotes standard $H^r(\Omega)$ norm (for r $\geq$ 0 and integer), where $H^r(\Omega)$ is the standard Sobolev space of order r,2 on $\Omega$. For r=0 this is $L_2$ norm. Similar to our previous consideration the superscript on $\|\cdot\|_r$ denotes the domain of integration other than $\Omega$. \\
The above result can be applied for the another interpolation error $\lambda_\eta$ similarly, i.e.
\begin{equation}
\|\lambda_\eta\|_s \leq \bar{C}(p,\Omega) h^{p+1-s} \|\eta\|_{p+1}
\end{equation}
Another results are strong stability conditions \cite{RefE},
\begin{equation}
\|c\|_2 \leq C_s \|q\|
\end{equation}
and similarly for the dual problem
\begin{equation}
\|\eta\|_2 \leq C_s' \|e\|
\end{equation}
Now looking at the RHS of the equation (16), we see that $a(e,\lambda_\eta)$ is the Galerkin part, which can be bounded by the result obtained in \cite{RefE}. But it is not possible to find the desired bound for the first term due to presence of $e_h$ . Therefore to get rid of $e_h$ we are doing the following steps.\\
Replacing $c_h$ by $c_h'+ e_h$ in SGS formulation (6) we have
\begin{equation}
\begin{split}
a_{SGS} (c_h'+e_h,d_h) &=l_{SGS} (d_h) \quad \forall d_h \in V_h\\
a_{SGS} (c_h',d_h)+ a_{SGS}(e_h,d_h) &=l_{SGS} (d_h) \quad \forall d_h \in V_h \\
a_{SGS} (c_h',d_h)+ a_{SGS}(e,d_h)- a_{SGS}(\lambda,d_h) & =l_{SGS} (d_h) \quad \forall d_h \in V_h\\
a_{SGS} (c_h',d_h)- a_{SGS}(\lambda,d_h) & =l_{SGS} (d_h) \quad \forall d_h \in V_h\\
\end{split}
\end{equation}
Subtracting (21) from (8), we will have
\begin{equation}
\begin{split}
a_{SGS} (e_h,d_h)+ a_{SGS}(\lambda,d_h) & =0 \quad \forall d_h \in V_h\\
a(e_h,d_h) + (-\mathcal{L^*} d_h,\tau \mathcal{L}e_h) + a(\lambda,d_h) + (-\mathcal{L^*} d_h,\tau \mathcal{L}\lambda) &= 0 \quad \forall d_h \in V_h\\
(\mathcal{L^*} d_h,\tau \mathcal{L}e_h) + (\mathcal{L^*} d_h,\tau \mathcal{L}\lambda) & = 0 \quad \forall d_h \in V_h\\
\end{split}
\end{equation}
Since $\eta_h$ is in $V_h$ therefore the above equation can be rewritten as:
\begin{equation}
(\mathcal{L^*} \eta_h,\tau \mathcal{L}e_h) + (\mathcal{L^*} \eta_h,\tau \mathcal{L}\lambda)=0 \\
\end{equation}
Applying the result obtained above into equation (16), we will remain with the Galerkin part only. As mentioned above, we now use the result from \cite{RefE} obtained to bound the Galerkin part using the standard interpolation and stability estimates. Therefore we finally have,
\begin{equation}
\begin{split}
\|e\|^2 & \leq \mid a(e,\lambda_\eta)\mid \\
        & \leq C h^{p+1} \|e\| \|c\|_{p+1}
\end{split}
\end{equation}
Hence,
\begin{equation}
\boxed{\|e\| \leq C h^{p+1} \|c\|_{p+1}}
\end{equation}
\begin{remark}
If we choose the finite element space as P1, the space containing continuous piecewise polynomial of order 1, a priori estimate using the stability estimate (19) will be as follows:
\begin{equation}
\begin{split}
\|e\| & \leq C h^{2} \|c\|_{2} \\
      & \leq C C_s h^2 \|q\|
\end{split}
\end{equation}
\end{remark}

\section{A posteriori error estimation in $L_2$ norm}
This estimation is carried out after the solution is known and it depends upon that computed solution, therefore it is computable. It ensures the order of convergence obtained in a priori estimate. Here we are going to find residual based a posteriori estimate. Residual denoted by $r_h$ is the term $(\mathcal{L}c_h - q)$. We will start with the same Nitsche trick i.e. considering the dual problem to proceed further. The steps will be slightly different from what have been done in previous section because we have to get rid of the exact solution c and bring residual into the picture. In this section again we have to deal with the Galerkin part, bound for which will be taken from \cite{RefE}.\\
Now let us again introduce the dual problem as follows: Find $\eta \in V$ such that 
\begin{equation}
a(d,\eta)= (d,e) \quad \forall d \in V
\end{equation}
According to the definition of $L_2$ norm,
\begin{equation}
\begin{split}
\|e\|^2 & = (e,e)\\
        & = a(e,\eta)\\
        & = a(c_h-c, \eta)\\
        & = a(c_h, \eta)- a(c, \eta)\\
        & = a(c_h, \eta) - l(\eta)\\ 
\end{split}
\end{equation}
Using weak form (4) in the last line by the fact that $\eta \in V$ we get rid of the exact solution.\\
The SGS stabilized form (8) can be rewritten as follows: Find $c_h \in V_h$ such that 
\begin{equation}
\begin{split}
a(c_h,d_h) + (-\mathcal{L}^*d_h, \tau ((\mathcal{L}c_h - q))_{\Omega'} & = l(d_h) \quad \forall d_h \in V_h \\
a(c_h,d_h) + (-\mathcal{L}^*d_h, \tau r_h)_{\Omega'} & = l(d_h) \quad \forall d_h \in V_h\\     
\end{split}
\end{equation}
Let $\eta_h$ be the interpolant of $\eta$ in $V_h$. Hence the above equation can be rewritten as:
\begin{equation}
a(c_h,\eta_h) + (-\mathcal{L}^*\eta_h, \tau r_h)_{\Omega'} - l(\eta_h) = 0
\end{equation}
Here we follow the same notation $\lambda_{\eta}$ for the interpolation error $(\eta_h - \eta)$.\\ 
Before proceeding further let us introduce the following appropriate interpolation estimate \cite{RefH} for each element $\Omega_k$, k=1,2,...,$n_{el}$
\begin{equation}
\begin{split}
\|\lambda_{\eta}\|_{s,\Omega_k} & \leq \bar{C}(\Omega_k) \hspace{1mm} h_k^{2-s} \hspace{1mm} \|\eta\|_2 \\
& \leq \bar{C}(\Omega_k) \hspace{1mm} C_s' \hspace{1mm} h_k^{2-s} \hspace{1mm} \|e\|
\end{split}
\end{equation}
where s is an integer and the second step comes from using the strong stability condition (20).\\
Now we are coming back in deducing the estimation by subtracting (30) from (28),
\begin{equation}
\begin{split}
\|e\|^2 & = a(c_h,\eta -\eta_h) - (-\mathcal{L}^*\eta_h, \tau r_h)_{\Omega'} - l(\eta -\eta_h)\\
& = (\mathcal{L}^*\eta_h, \tau r_h)_{\Omega'} + l(\eta_h -\eta) - a(c_h,\eta_h -\eta)\\
& = (\mathcal{L}^*\eta_h, \tau r_h)_{\Omega'} + l(\lambda_\eta) - a(c_h,\lambda_\eta)\\  
& \leq \mid (\mathcal{L}^*\eta_h, \tau r_h)_{\Omega'} \mid + \mid a(c_h,\lambda_\eta)- l(\lambda_\eta) \mid \\
\end{split}
\end{equation}
Here the second term is the Galerkin part. \cite{RefE} gives its bound using the above interpolation estimate as follows:
\begin{equation}
\mid a(c_h,\lambda_\eta)- l(\lambda_\eta) \mid \leq (\sum_{k=1}^{n_{el}} \bar{C}(\Omega_k) \hspace{1mm} C_s' \hspace{1mm} h_k^2 \hspace{1mm} \|r_h\|_{\Omega_k} )\hspace{1mm} \|e\|
\end{equation}
Now our job is to find bound for the first term.
\begin{equation}
\begin{split}
 \mid (\mathcal{L}^*\eta_h, \tau r_h)_{\Omega'} \mid & =  \mid (\mathcal{L}^*\lambda_\eta, \tau r_h)_{\Omega'} + (\mathcal{L}^*\eta, \tau r_h)_{\Omega'} \mid \\
 & \leq \mid (\mathcal{L}^*\lambda_\eta, \tau r_h)_{\Omega'} \mid + \mid (e , \tau r_h)_{\Omega'} \mid \\
 & = \mid \tau \mid \hspace{1mm} \mid \sum_{k=1}^{n_{el}}(\mathcal{L}^*\lambda_\eta, r_h)_{\Omega_k} \mid + \mid \tau \mid \hspace{1mm} \mid \sum_{k=1}^{n_{el}} (e , r_h)_{\Omega_k} \mid \\
 & \leq \mid \tau \mid (\sum_{k=1}^{n_{el}} \| \mathcal{L}^*\lambda_\eta \|_{\Omega_k} \|r_h\|_{\Omega_k} + \sum_{k=1}^{n_{el}} \|e\|_{\Omega_k} \|r_h\|_{\Omega_k}) \\
& \leq \mid \tau \mid (\sum_{k=1}^{n_{el}} \| \mathcal{L}^*\lambda_\eta \|_{\Omega_k} \|r_h\|_{\Omega_k} + \sum_{k=1}^{n_{el}}  \|r_h\|_{\Omega_k} \|e\|) 
\end{split}
\end{equation}
Now we are only remained with estimating $\|\mathcal{L}^*\lambda_\eta\|_{\Omega_k}$ as follows:
\begin{equation}
\begin{split}
\|\mathcal{L}^*\lambda_\eta\|_{\Omega_k} & = \| -\bigtriangledown \cdot \widetilde{\bigtriangledown}\lambda_\eta - \overline{u} \cdot \bigtriangledown \lambda_\eta + \mu \lambda_\eta \|_{\Omega_k} \\
& = \|-D_1 \frac{\partial^2 \lambda_\eta}{\partial x^2}- D_2 \frac{\partial^2 \lambda_\eta}{\partial y^2} +(u_1- \frac{\partial D_x}{\partial x}) \frac{\partial \lambda_\eta}{\partial x} +(u_2- \frac{\partial D_y}{\partial y}) \frac{\partial \lambda_\eta}{\partial y} + \mu \lambda_\eta \|_{\Omega_k} \\
& \leq D \hspace{1mm} (\|\frac{\partial^2 \lambda_\eta}{\partial x^2}\|_{\Omega_k} + \| \frac{\partial^2 \lambda_\eta}{\partial y^2} \|_{\Omega_k} + \|\frac{\partial \lambda_\eta}{\partial x}\|_{\Omega_k} + \|\frac{\partial \lambda_\eta}{\partial y}\|_{\Omega_k} + \|\lambda_\eta\|_{\Omega_k} )\\
& = D \hspace{1mm} \|\lambda_\eta\|_{2,{\Omega_k}} \quad (by \hspace{1mm} standard \hspace{1mm} definition \hspace{1mm} of \hspace{1mm} H^2 \hspace{1mm} norm) \\
& \leq D \hspace{1mm} \bar{C}(\Omega_k) \hspace{1mm} C_s' \hspace{1mm} \|e\|
\end{split}
\end{equation}
where D is modulus of the maximum taken over all the coefficients present in that expression and we have obtained the last line using interpolation estimate (31).\\
Putting the result obtained in (35) into (34), we will have,
\begin{equation}
\begin{split}
\mid (\mathcal{L}^*\eta_h, \tau r_h)_{\Omega'} \mid & \leq (\sum_{k=1}^{n_{el}}( D \hspace{1mm} \bar{C}(\Omega_k) \hspace{1mm} C_s' + 1) \hspace{1mm} \|r_h\|_{\Omega_k})\hspace{1mm} \mid \tau \mid \|e\| \\
& = \sum_{k=1}^{n_{el}} \bar{D}_k \mid \tau \mid \|r_h\|_{\Omega_k} \hspace{1mm} \|e\| \\
\end{split}
\end{equation}
where $\bar{D}_k = ( D \hspace{1mm} \bar{C}(\Omega_k) \hspace{1mm} C_s' + 1)$ \\
Combining both the results obtained in (33) and (36) and substituting them into equation (32), we will remain with the following:
\begin{equation}
\|e\|^2 \leq (\sum_{k=1}^{n_{el}} ( \bar{C}(\Omega_k) \hspace{1mm} C_s' \hspace{1mm} h_k^2 \hspace{1mm}+ \bar{D}_k \mid \tau \mid ) \hspace{1mm}\|r_h\|_{\Omega_k})\hspace{1mm} \|e\| 
\end{equation}
Hence,
\begin{equation}
\boxed {\|e\| \leq (\sum_{k=1}^{n_{el}} ( \bar{C}(\Omega_k) \hspace{1mm} C_s' \hspace{1mm} h_k^2 \hspace{1mm}+ \bar{D}_k \mid \tau \mid ) \hspace{1mm}\|r_h\|_{\Omega_k})}
\end{equation}
\begin{remark}
The above a posteriori estimate is computable with an appropriate estimate of the constants $ \bar{C}(\Omega_k), C_s', \bar{D}_k $ and the unknown stabilization parameter $\tau$. An estimation of $\tau$ is derived in the next section.
\end{remark}
\subsection{Estimating $\tau$}
In \cite{RefB} and \cite{RefC}, Codina gives an approach based on discrete maximum principle (DMP) to find conditions on $\tau$ for scalar 1D case. We are going to apply the same idea to reach to a condition on $\tau$ for general 2D case and then use that for 2D ADR equation with spatially variable coefficients. The idea has been taken up from the following implications:\\
The finite element matrix A in non-negative $\Longrightarrow$ the DMP holds $\Longrightarrow$ stability and convergence of the finite element solution. \vspace{2 mm}\\
In both the references mentioned in the above paragraph, a condition on $\tau $ has been obtained by showing the matrix of the final algebraic system to be non-negative. Before introducing the definition of non-negativity of matrix A, let us mention some useful notions. We already have $n_{pt}$ as the total number of nodes of finite element mesh. Let $n_{free}$ be the interior nodes. The finite element discretization of the problem leads to (7), where the values $c_a$ for a= $n_{free}$+1,...,$n_{pt}$ i.e. for the boundary nodes, are known from the Dirichlet boundary conditions. Now the matrix A have dimensions $n_{free}$ $\times$ $n_{pt}$ and the column vector on the RHS will have $n_{free}$ components.\\
\begin{defn}
Matrix A= $(a_{ij})_{n_{free} \times n_{pt}}$ is of non-negative type if the following conditions hold:\\
(i) $a_{ij} \leq 0$ $\forall $ $i=1,2,...,n_{free}$ , $j=1,2,...,n_{pt}$ ,$i \neq j$ \\
(j) $\sum_{j=1}^{n_{pt}} a_{ij} \geq 0$ for $i=1,2,...,n_{free}$
\end{defn}
To find a condition on $\tau$ requiring the finite element matrix to be non-negative we first take  a sup norm of the spatially variable coefficient functions because if the condition holds for the supremum of the coefficients, it will hold for every possible values of the coefficients as well as we can deal with the variational form of constant coefficients. The numerical experiments will explain this better later.\\
Let D and U be the sup norm values taken over the diffusion coefficients and velocity components respectively on $\Omega$.\\ 
Hence the variational form (4) will look like as follows: Find c $\in$ V such that 
\begin{equation}
B(c,d) = F(d) \forall d \in V
\end{equation} 
where \\

\hspace{0.7cm} B(c,d)=$D \int_{\Omega}(\frac{\partial c}{\partial x}\frac{\partial d}{\partial x} +\frac{\partial c}{\partial y} \frac{\partial d}{\partial y}) + U \int_{\Omega}d (\frac{\partial c}{\partial x} + \frac{\partial c}{\partial y})  + \mu \int_{\Omega} c d  $\\ 

\hspace{1cm} l(d)=$\int_{\Omega} q d = (d,q) $\\
Since the assembly operator is summation of the element contributions, therefore it is sufficient to check the element matrices to be non-negative. Hence we will now try to find the condition for one element $\Omega_e$.\vspace{1 mm}\\
Let us discretize the domain into triangular elements and $\Lambda_e$ be the area of each element $\Omega_e$ for e=1,2,...,$n_{el}$. We have taken continuous piecewise linear polynomials as the shape functions for each element. Letting $\{N_1^e, N_2^e, N_3^e \}$ be the set of basis functions for each element $\Omega_e$, the element matrices will be as follows:\\
\[\begin{bmatrix}
B^e(N_1^e,N_1^e)&B^e(N_1^e,N_2^e)&B^e(N_1^e,N_3^e)\\
B^e(N_2^e,N_1^e)&B^e(N_2^e,N_2^e)&B^e(N_2^e,N_3^e)\\
B^e(N_3^e,N_1^e)&B^e(N_3^e,N_2^e)&B^e(N_3^e,N_3^e)\\
\end{bmatrix}\]
where $B^e(\cdot,\cdot)$ is bilinear form defined on each element $\Omega_e$\\
$N_1^e (x,y)= \frac{1}{2 \Lambda_e} \{ (x_2^e y_3^e - x_3^e y_2^e) + (y_2^e-y_3^e) x + (x_3^e- x_2^e) y \}$\\
$N_2^e (x,y)= \frac{1}{2 \Lambda_e} \{ (x_3^e y_1^e - x_1^e y_3^e) + (y_3^e-y_1^e) x + (x_1^e- x_3^e) y \}$\\
$N_3^e (x,y)= \frac{1}{2 \Lambda_e} \{ (x_1^e y_2^e - x_2^e y_1^e) + (y_1^e-y_2^e) x + (x_2^e- x_1^e) y \}$ \vspace{1 mm}\\
where $(x_1^e, y_1^e), (x_2^e, y_2^e)$ and $(x_3^e, y_3^e)$ are three vertices of element $\Omega_e$.\\

Introducing the notations $\beta_i^e = (y_j^e-y_k^e)$ and $\gamma^e = -(x_j^e-x_k^e)$, the orientation of vertices can be considered such a way in anti-clock wise direction that the following form will hold for each element.\\

$B^e(N_i^e,N_j^e) = - \frac{D}{4 \Lambda_e} K_{ij}^e \pm \frac{U}{6} E_{i}^e + \frac{\mu \Lambda_e}{12}$ for i,j=1,2,3. \vspace{2 mm}\\
here $K_{ij}^e = (\beta_i^e \beta_j^e + \gamma_i^e \gamma_j^e)$ and $E_{i}^e = (\beta_i^e + \gamma_i^e)$ both are always positive constants.\\
According to the definition (1), non-negativity of each element matrix is equivalent to
\begin{equation}
- \frac{D}{4 \Lambda_e} K_{ij}^e \pm \frac{U}{6} E_{i}^e + \frac{\mu \Lambda_e}{12} \leq 0 \quad \forall i,j =1,2,3 \quad and \quad i \neq j
\end{equation}
When SGS method is applied, the finite element matrices will have the same form (40), but with modified coefficients such as:\\
$\bar{D} = D + \tau U^2 \quad \bar{U} = U - 2 \mu \tau U \quad \bar{\mu} = \mu - \tau \mu^2 $ \\
Hence the above condition (40) will be modified as follows:
\begin{equation}
(- \frac{D}{4 \Lambda_e} K_{ij}^e \pm \frac{U}{6} E_{i}^e + \frac{\mu \Lambda_e}{12}) - \tau (- \frac{U^2}{4 \Lambda_e} K_{ij}^e \pm \frac{2\mu U}{6} E_{i}^e + \frac{\mu^2 \Lambda_e}{12}) \leq 0 \quad \forall i,j=1,2,3 \quad and \quad
 i \neq j.
\end{equation}
We will reach at a condition on $\tau$ by considering two simplified cases following \cite{RefC}, as follows:\\
\textbf{Case 1}\\
Let us consider the reaction term $\mu$ = 0. Then (41) becomes,
\begin{equation}
\begin{split}
\tau & \geq (-\frac{D}{U^2} \pm \frac{2 \Lambda_e E_{i}^e}{3 U K_{ij}^e})\\
& \geq (\frac{2 \Lambda_e E_{i}^e}{3 U K_{ij}^e} - \frac{D}{U^2})\\
& \geq (\frac{2 h_e}{3 U } - \frac{D}{U^2})\\
& = \frac{2 h}{3 U } (1- \frac{1}{Q}) 
\end{split}
\end{equation}
where we take Q = $\frac{2Uh_e}{3D}$ and h satisfies the condition $\frac{\Lambda_e E_{i}^e}{K_{ij}^e} \geq h_e$\\
\textbf{Case 2}\\
This time let us consider the advection velocity U to be 0. Then the condition (41) reduces to,
\begin{equation}
\begin{split}
\tau & \geq (\frac{1}{\mu} - \frac{3D K_{ij}^e}{\mu^2 \Lambda_e^2})\\
& \geq (\frac{1}{\mu} - \frac{3D K_{ij}^e}{\mu^2 h_e^2})\\
& = \frac{1}{\mu} (1-\frac{3}{R})
\end{split}
\end{equation}
where $\frac{\mu h^2}{D}$ is taken as R along with a similar consideration as above, $\frac{\Lambda_e^2}{K_{ij}^e} \geq h_e^2$\\
Since, $\quad$
$\frac{2 h_e}{3 U } (1- \frac{1}{Q}) \leq \frac{2 h_e}{3 U } (\frac{Q}{1+Q})= (\frac{9D}{4h_e^2} + \frac{3U}{2h_e})^{-1} \leq (\frac{9D}{4h^2} + \frac{3U}{2h})^{-1} $\\
And, $\quad$ $\frac{1}{\mu} (1-\frac{3}{R}) \leq \frac{1}{\mu} (\frac{4 R}{4R+9}) = (\frac{9D}{4h_e^2} + \mu)^{-1} \leq (\frac{9D}{4h^2} + \mu)^{-1} $ \\
Now following the approach in \cite{RefC}, $\tau$ can be chosen as,
\begin{equation}
\boxed{\tau = (\frac{9D}{4h^2} + \frac{3U}{2h} + \mu )^{-1}}
\end{equation}
\begin{remark}
From the above expression of $\tau$ we can conclude that $\tau$ is of order of $h^2$.
\end{remark}
\section{Numerical Experiment}
As mentioned earlier, Galerkin finite element method fails when diffusion coefficients are small compared to either advection or reaction coefficients or both. In this section we will prove this fact as well as verify the theoretically established result through numerical data obtained by solving a test problem with spatially variable coefficients of hydrological significance. Here we have chosen two cases. Whereas the first case will analyse the comparison between the Galerkin and the stabilized SGS method for diffusion dominated flow, the second one will do the same for advection- reaction dominated flow.\vspace{2mm}\\
As a test problem we have considered a source of pollutant dispersed into an incompressible fluid in a simple bounded square domain $\Omega = (0,1) \times (0,1)$. This phenomena respects ADR equation with spatially variable coefficients along with homogeneous Dirichlet boundary condition. The concepts of choosing  the coefficients are taken from \cite{RefD}, which is motivated by the fact of real-life application from hydrology. In that study, velocity components are taken spatial linear non-homogeneous functions, whereas the diffusion coefficients are considered proportional to the square of respective velocity components.\\
\subsection{First case: Diffusion dominated flow}
In this case we will see that the Galerkin finite element method and the stabilized SGS method perform equally well in establishing convergence. Let us introduce the coefficients we have worked with as follows:\\
$u_1= 0.5(1+0.02x)$, $u_2 = -0.5(1+0.02y)$, $D_1= (1+0.02x)^2$, $D_2= 0.1(1+0.02y)^2$ and $\mu= 0.01$\\
Then D=1.0404 and U= 0.51.\\
The exact solution is c= sin(xy(x-1)(y-1)) on $\Omega$ along with c vanishing on the boundary $\partial\Omega$. Working with the same expression for $\tau$ as deduced in the previous section, we have obtained results using \textit{freefem++}.\\
\begin{remark}
The error and order of convergence obtained for the Galerkin and the SGS method for diffusion dominated flow have been shown in table 1 and table 2 respectively.
\end{remark}

\begin{remark}
The results of the table 2 shows that the order of convergence for the stabilized SGS method is 2 which also justifies the deduced expression the stabilization parameter $\tau$ numerically.
\end{remark}

\begin{table}
   \centering
    \begin{tabular}{||c c c||}
    \hline 
    Mesh size & Error in $L_2$ norm &  Order of convergence \\
    \hline \hline
      10      &  0.000254151        &                       \\
      \hline
      20      &  6.28651 $e^{-5}$   & 2.01536  \\
      \hline
      40      &  1.26455 $e^{-5}$   & 2.31363  \\
      \hline
      80      &  3.01483 $e^{-6}$   & 2.06848  \\
      \hline
      160     &  8.30888 $e^{-7}$   & 1.85935  \\
      \hline
      320     &  1.72972 $e^{-7}$   & 2.26412 \\
      \hline
    \end{tabular}
     \caption{Error and Order of convergence obtained in $L_2$ norm by Galerkin method for diffusion dominated flow}
   % \label{table:draglift1}
\end{table}

\begin{table}[]
    \centering
    \begin{tabular}{||c c c||}
    \hline 
    Mesh size & Error in $L_2$ norm &  Order of convergence \\
    \hline \hline
      10      &   0.000405122        &                       \\
      \hline
      20      &  0.000105007        & 1.94787 \\
      \hline
      40      &  2.40648 $e^{-5}$   & 2.12549  \\
      \hline
      80      &  5.89785 $e^{-6}$   & 2.02866  \\
      \hline
      160     &  1.55175 $e^{-6}$   & 1.92629  \\
      \hline
      320     &  3.49761 $e^{-7}$   & 2.14945 \\
      \hline
    \end{tabular}
\caption{Error and Order of convergence obtained in $L_2$ norm by the SGS method for diffusion dominated flow}
    %\label{table:draglift2}
\end{table} 
\begin{table}
   \centering
    \begin{tabular}{||c c c||}
    \hline 
    Mesh size & Error in $L_2$ norm &  Order of convergence \\
    \hline \hline
      10      &  0.000541066        &                       \\
      \hline
      20      &  0.000375131        & 0.528409  \\
      \hline
      40      &  9.51933 $e^{-5}$   & 1.97846 \\
      \hline
      80      &  3.399958 $e^{-5}$   & 1.4855  \\
      \hline
      160     &  1.25028 $e^{-5}$   & 1.4431  \\
      \hline
      320     &  4.00531 $e^{-6}$   & 1.64227 \\
      \hline
    \end{tabular}
     \caption{Error and Order of convergence obtained in $L_2$ norm by Galerkin method for flow having small diffusion coefficients}
    %\label{table:draglift3}
\end{table}

\begin{table}[]
    \centering
    \begin{tabular}{||c c c||}
    \hline 
    Mesh size & Error in $L_2$ norm &  Order of convergence \\
    \hline \hline
      10      &   0.000280761        &                       \\
      \hline
      20      &  6.90662 $e^{-5}$   & 2.02329 \\
      \hline
      40      &  1.73534 $e^{-5}$   & 1.99276  \\
      \hline
      80      &  4.64904 $e^{-6}$   & 1.90022  \\
      \hline
      160     &  1.37837 $e^{-6}$   & 1.85397  \\
      \hline
      320     &  3.50102 $e^{-7}$   & 1.97711 \\
      \hline
    \end{tabular}
\caption{Error and Order of convergence obtained in $L_2$ norm by the SGS method for flow having small diffusion coefficients}
    %\label{table:draglift4}
\end{table}

\subsection{Second case: Advection-reaction dominated flow}
We have chosen diffusion very small compared to advection and reaction coefficients. The coefficients following the same pattern as above are,\\
$u_1= (1+0.02x)$, $u_2 = -(1+0.02y)$, $D_1= 10^{-7} (1+0.02x)^2$, $D_2= 10^{-8}(1+0.02y)^2$ and $\mu= 1$\\
Then D= 1.0404 $\times 10^{-7}$ and U=1.02\\
The exact solution will be the same.\\
\begin{remark}
Table 3 and table 4 represent the error and order of convergence obtained for Galerkin and the SGS method  respectively for very small diffusion coefficients. Table 3 clearly shows oscillations in order of convergence for Galerkin method whereas The order of convergence under the SGS method is converging to 2.
\end{remark}
%\section{Conclusion}

\end{document}